\magnification \magstep1
\hoffset .13in
\def\fpee{{{\bf F}_p}}
\def\dt{{*}}
\def\mod{|}
\def\longarrowright{\rightarrow}
\def\pra{\par}

\def\be{1}
\def\ca{2}
\def\ce{3}
\def\di{4}
\def\ev{5}
\def\hm{6}
\def\ka{7}
\def\ku{8}
\def\se{9}
\def\st{10}
\def\cupio{{\smile_{1,0}}}
\def\cupoi{{\smile_{0,1}}}
\def\pr{{\wedge}}
\def\proof{\noindent{\bf Proof.}\quad}
\def\remark{\noindent{\bf Remark.}\quad}
\def\qed{\qquad{$QED$}}

\centerline{\bf A differential in the Lyndon-Hochschild-Serre spectral 
sequence}
\centerline{by I. J. Leary\footnote{*}{supported by SERC post-doctoral 
fellowship number B90 RFH 8960}}
\centerline{School of Mathematical Sciences,}
\centerline{Queen Mary and Westfield College,}
\centerline{Mile End Road,}
\centerline{London}
\centerline{E1 4NS}
\par 
\noindent

\beginsection Abstract.

We consider the Lyndon-Hochschild-Serre spectral sequence with $\fpee$
coefficients for a central extension with kernel cyclic of order a
power of $p$ and arbitrary discrete quotient group.  For this spectral
sequence $d_2$ and $d_3$ are known, and we give a description
for $d_4$.  Using this result we deduce a similar formula for the Serre 
spectral sequence for a fibration with fibre $BC_{p^n}$.  
The differential from odd rows to even rows involves a Massey triple
product, so we describe the calculation of such products in the
cohomology of a finite abelian group.  As an example we determine the
Poincar\'e series for the mod-3 cohomology of various 3-groups.  
\par
\beginsection Introduction. 

Let $G$ be a discrete group, let $C$ be a cyclic group of order $p^m$
for some prime $p$ and $m\geq 1$, and let $E$ be a central extension
with kernel $C$ and quotient $G$.  Such an extension is described by
its extension class, which may be any element of $H^2(G;C)$.  Since
$C$ is central, the $E_2$-page of the Lyndon-Hochschild-Serre spectral 
sequence with $\fpee$ coefficients for this extension is isomorphic to 
the graded tensor product of $H^*(G)$ and $H^*(C)$.  (Note that throughout
this paper $H^*(.)$ will denote mod-$p$ cohomology.)  There is an isomorphism
$$H^*(C)\cong \fpee[u,t]/I,$$
where $u$ has degree one, $t$ has degree two, and $I$ is generated by 
$u^2$ except when $p^m=2$ in which case $I$ is generated by $u^2-t$.  
In particular this implies that each row of the $E_2$-page of 
the spectral sequence is isomorphic to a single copy of $H^*(G)$.  
\par
The first two differentials in the spectral sequence have been known
since the calculations of the cohomology of Eilenberg-MacLane spaces 
by Cartan and Serre [\ca], [\se].  Modulo choice of various
isomorphisms, $d_2(u)$ is the mod-$p$ reduction of the extension class,
$d_2(t)$ is zero, and $d_3(t)$ is the image of the extension class
under the Bockstein for the following coefficient sequence.
$$0\longarrowright \fpee\longarrowright {{\bf Z}/{p^{m+1}}}
\longarrowright {{\bf Z}/{p^m}}\cong C\longarrowright 0 $$
We shall describe the next differential $d_4$, in Theorem~3, but first 
we require two lemmata concerning the spectral sequence.  
In Corollary~4 we deduce a similar result for the Serre spectral
sequence for a fibration with fibre $BC_{p^m}$.

\beginsection The spectral sequence.

We set up the Lyndon-Hochschild-Serre spectral sequence as described
in Evens' book [\ev].  Thus we let $P_\dt$ be the standard (or bar)
resolution for $\fpee$ over $\fpee G$, and $Q_\dt $ be the standard 
resolution for $\fpee$ over $\fpee E$.  Now we define 
$$E_0^{i,j}= {\rm Hom}_E(P_i\otimes Q_j,\fpee),$$
where $E$ acts diagonally on $P_\dt\otimes Q_\dt$, $d_1$ is the map
induced by the differential on $P_\dt$, and $d_0$ is the map induced 
by the differential on $Q_\dt$ except that a sign must be introduced
(in fact on $E_0^{i,j}$ $d_0$ is the adjoint of $(-1)^i1\otimes d^Q$).
Using the Kunneth theorem it may be checked that the singly graded 
total complex $(P_\dt\otimes Q_\dt)_\dt$ with differential $d_0+d_1$
is an $E$-free resolution for $\fpee$.  
\par
To introduce products into the spectral sequence we need a diagonal 
approximation for $(P_\dt\otimes Q_\dt)_\dt$, that is a chain map from 
$(P_\dt\otimes Q_\dt)_\dt$ to $(P_\dt\otimes Q_\dt\otimes P_\dt\otimes
Q_\dt)_\dt$ inducing the identity map on homology.  For any two chain 
complexes $A_\dt $ and $B_\dt$, let $\tau$ denote the chain
isomorphism from $A_\dt\otimes B_\dt$ to $B_\dt \otimes A_\dt$ given
by $\tau(a\otimes b)= (-1)^{\mod a\mod\mod b\mod}$.  If $\Delta_0^P$ and
$\Delta_0^Q$ are diagonal approximations for $P_\dt$ and $Q_\dt$ respectively,
then $(1\otimes \tau\otimes 1)(\Delta_0^P\otimes\Delta_0^Q)$ is a 
diagonal approximation for $P_\dt\otimes Q_\dt$, which is
coassociative provided that $\Delta_0^P$ and $\Delta_0^Q$ are.  Since
it is possible to choose $\Delta_0^P$ and $\Delta_0^Q$ to be coassociative
using the Alexander-Whitney formula (given on the usual $\fpee$ basis
of $P_n$ by 
$$\Delta_0^P(g_0,\ldots,g_n)=\sum^n_{i=0}(g_0,\ldots,g_i)\otimes(g_i,
\ldots,g_n)$$
and similarly for $\Delta_0^Q$) we do so.  We thus obtain an
associative product (which we shall call the cup product) 
on the $E_0$-page of the spectral sequence given by 
$$\phi\smile\theta(x)=\mu(\phi\otimes\theta)(1\otimes\tau\otimes 1)
(\Delta_0^P\otimes\Delta_0^Q)(x)$$
for all $\phi$ and $\theta$ in $E_0$ and all $x\in P_\dt\otimes
Q_\dt$, where $\mu$ is the product map from $\fpee\otimes\fpee$ to
$\fpee$.  The only reason why we insisted that $P_\dt$ and $Q_\dt$
should be the standard resolution was to ensure that $E_0$ could be
endowed with an associative product.  The differentials $d_0$ and
$d_1$ are derivations for this product.  
\par
We also need to introduce analogues within the spectral sequence of 
the cup-1 product [\st].  The map $\tau \Delta_0^P$ is also a diagonal 
approximation for $P_\dt$, so is chain homotopic to $\Delta_0^P$.  Let 
$\Delta_1^P$ be a chain homotopy from $\Delta_0^P$ to
$\tau\Delta_0^P$.  For example,  $\Delta_1^P$ could be the map given
by Steenrod [\st]
$$\Delta_1^P(g_0,\ldots,g_n)=\sum_{0\leq i<j\leq n}(-1)^{f(i,j,n)}
(g_0,\ldots,g_i,g_j,\ldots,g_n)\otimes(g_i,\ldots,g_j),$$
where $f(i,j,n)=n + (n-i-1)(j-i-1)$.  We define $\Delta_1^Q$
similarly.  We now consider three different diagonal approximations for
$P_\dt\otimes Q_\dt$, which are defined as follows.
\def\sw{(1\otimes\tau\otimes 1)}
\def\do{(\Delta_0^P\otimes\Delta_0^Q)}
$$\sw\do \qquad \sw(\tau\otimes1\otimes1)\do\qquad \sw(\tau\otimes\tau)\do$$
As above, each of these induces a product on the $E_0$-page of the
spectral sequence.  By definition the first induces the cup product,
and it is easily checked that the third induces the product 
$$ (\phi,\theta)\mapsto(-1)^{|\phi||\theta|}\theta\smile\phi,$$
where $|.|$ refers to total degree.  The second induces another
product which we shall refer to as 
$$(\phi,\theta)\mapsto\phi\pr\theta.$$
It is easily checked that $\sw(\Delta_1^P\otimes\Delta_0^Q)$ is a
chain homotopy between the first two diagonal approximations, and that 
the map defined to be $(-1)^i\sw(\tau\otimes1\otimes1)(\Delta_0^P
\otimes\Delta_1^Q)$ on $P_i\otimes Q_j$ is a chain homotopy between
the second two.  We may use these chain homotopies to define cup-1
products of bidegrees $(-1,0)$ and $(0,-1)$ respectively, whose properties
are summarised in the following lemma.  
\par
\proclaim Lemma 1.  For $\phi, \theta\in E_0$, and $x\in P_i\otimes
Q_j$, define products of $\phi$ and $\theta$ as below.  
$$\eqalign{\phi\cupio\theta(x)&=\mu(\phi\otimes\theta)\sw
(\Delta_1^P\otimes\Delta_0^Q)(x)\cr
\phi\cupoi\theta(x)&=(-1)^i\mu(\phi\otimes\theta)\sw(\tau\otimes1\otimes1)
(\Delta_0^P\otimes\Delta_1^Q)(x)}$$
Then $\cupio$ and $\cupoi$ satisfy the following coboundary formulae.  
$$\eqalign{d_0(\phi\cupio\theta)&=
-(d_0\phi)\cupio\theta-(-1)^{|\phi|}\phi\cupio(d_0\theta)\cr
d_1(\phi\cupio\theta)&=
-(d_1\phi)\cupio\theta-(-1)^{|\phi|}\phi\cupio(d_1\theta)+\phi\smile\theta
-\phi\pr\theta\cr
d_0(\phi\cupoi\theta)&=
-(d_0\phi)\cupoi\theta-(-1)^{|\phi|}\phi\cupoi(d_0\theta)
+\phi\pr\theta-(-1)^{|\phi||\theta|}\theta\smile\phi\cr
d_1(\phi\cupoi\theta)&=
-(d_1\phi)\cupoi\theta-(-1)^{|\phi|}\phi\cupoi(d_1\theta).}$$
\par
\noindent
\proof
Left as an exercise for the reader.  \qed
\par
\beginsection Differentials. 

In any spectral sequence, if $\phi \in E_0^{i,j}$ is to survive until
the $E_2$-page, $d_0(\phi)$ must be zero, and also $d_1(\phi)$
must represent zero on the $E_1$-page, or equivalently there must be 
$\phi_1\in E_0^{i+1,j-1}$ such that $d_0(\phi_1)=d_1(\phi)$.  In this 
case $d_2(\phi)$ is represented by $d_1(\phi_1)$.  Similarly, $\phi$
can only survive until the $E_3$-page if there is 
$\phi_2\in E_0^{i+2,j-2}$ such that $d_1(\phi_1)=d_0(\phi_2)$, in
which case $d_3(\phi)$ is represented by $d_1(\phi_2)$.  In this way 
the higher differentials in the spectral sequence may be computed. 
\par
Returning to the spectral sequence we are studying,  let $u$ and $t$
be cochains in $E_0^{0,1}$ and $E_0^{0,2}$ respectively, representing 
the elements of $E_2^{0,*}$ with the same names (as used in the introduction).
The known description of $d_2$ implies that $d_0(u)=0$, and there is
an element $\theta$ in $E_0^{1,0}$ such that $d_1(u)=d_0(\theta)$, and 
$d_1(\theta) = \xi$ is a cocycle representing $d_2(u)$.  Similarly,
there are cochains $\eta_1\in E_0^{1,1}$ and $\eta_2\in E_0^{2,0}$
such that $d_0(t)=0$, $d_1(t)=d_0(\eta_1)$,
$d_1(\eta_1)=d_0(\eta_2)$, and $d_1(\eta_2)=\xi'$, where $\xi'$ is a
cocycle representing $d_3(t)$.  Note that $\xi$ and $\xi'$ are
universal cycles.  The situation may be summarised in the following diagram.  
$$\matrix{0&& && && 0 \cr
\uparrow&& && && \uparrow \cr
u&\rightarrow&\bullet&& && t&\rightarrow&\bullet \cr
&&\uparrow&& && && \uparrow \cr
&&\theta &\rightarrow &\xi&& && \eta_1&\rightarrow&\bullet\cr
&& && && && &&  \uparrow\cr
&& && && && &&  \eta_2&\rightarrow&\xi'\cr}$$
From now on we fix this notation, so that $u$, $t$, $\theta$, $\xi$, 
$\eta_1$, $\eta_2$ and $\xi'$ will always refer to these cochains.  
Recall from the introduction that the classes represented by $\xi$ and
$\xi'$ in $H^*(G)$ may be easily described in terms of the extension
class of the group $E$ in $H^2(G;C)$.  We are now able to state our 
second lemma, which describes a sequence of cochains related to $t^n$.
\par
\proclaim Lemma 2.  For all $n \geq 1$ there exist cochains $\eta_1(n),
\ldots,\eta_4(n)$, with $\eta_i(n)\in E_0^{i,2n-i}$ such that 
$$\eqalign{d_1(t^n)&=d_0(\eta_1(n))\cr
d_1(\eta_2(n))&=nt^{n-1}\xi'+d_0(\eta_3(n))}\qquad
\eqalign{d_1(\eta_1(n))&=d_0(\eta_2(n))\cr
d_1(\eta_3(n))&=n\eta_1(n-1)\xi'+d_0(\eta_4(n)).}$$
\par\noindent
\remark The existence of $\eta_1(n)$, $\eta_2(n)$ and $\eta_3(n)$ and
the relations between them follow immediately from the fact that 
$d_3(t^n)$ is represented by $nt^{n-1}\xi'$.  Thus the content of this
lemma is the existence of $\eta_4(n)$ and the equation for $d_1(\eta_3(n))$.
\par\noindent
\proof  Define $\eta_i(0)=\eta_3(1)=\eta_4(1)=0$, $\eta_1(1)=\eta_1$, and
$\eta_2(1)=\eta_2$, where $\eta_1$ and $\eta_2$ are as above.  Then define
$\eta_i(n)$ inductively by the following formulae.
$$\eqalign{\eta_1(n)=&\eta_1(n-1)t+t^{n-1}\eta_1\cr
\eta_2(n)=&\eta_2(n-1)t+\eta_1(n-1)\eta_1+t^{n-1}\eta_2 
-(n-1)t^{n-2}(\xi'\cupio t)\cr
\eta_3(n)=&\eta_3(n-1)t+\eta_2(n-1)\eta_1+\eta_1(n-1)\eta_2\cr
&-(n-1)(\eta_1(n-2)(\xi'\cupio t)+t^{n-2}(\xi'\cupio\eta_1))
+(n-1)t^{n-2}(\eta'\cupoi t)\cr
\eta_4(n)=&\eta_4(n-1)t+\eta_3(n-1)\eta_1+\eta_2(n-1)\eta_2\cr
&-(n-1)(\eta_2(n-2)(\xi'\cupio t)+\eta_1(n-2)(\xi'\cupio \eta_1)
+t^{n-2}(\xi'\cupio \eta_2)) \cr
&+(n-1)(\eta_1(n-2)(\xi'\cupoi t)+t^{n-2}(\xi'\cupoi \eta_1))}$$
The equations given in the statement hold for $n=1$, and may be
verified for all $n$ by induction.  \qed
\par
We are now ready to consider $d_4$.  
First we decide which elements of $E_2$ will survive until $E_4$.  Let
$\chi$ be an element of $H^*(G)$.  Since $d_2(t^nu\chi)=t^n\xi\chi$
we deduce that $t^nu\chi$ survives until $E_3$ iff $\xi\chi=0$ in
$H^*(G)$.  In this case $d_3(t^nu\chi)=-nt^{n-1}u\xi'\chi$, and since  
$d_2$ had trivial image in $E_2^{*,2n-1}$ it follows that $t^nu\chi$
survives until $E_4$ iff $\xi\chi=0$ in $H^*(G)$ and either $p$ 
divides $n$ or $\xi'\chi=0$ in $H^*(G)$.  If $p$ divides $n$ then 
since $u\chi$ is a
universal cycle and $t^p$ survives until $E_{2p+1}$ (by the Serre
transgression theorem) then $d_4(t^nu\chi)=0$.  Thus we only need
determine $d_4(t^nu\chi)$ for $\chi$ such that $\chi\xi$ and
$\chi\xi'$ are zero.  Similarly, for $t^n\chi$ to survive, we must
have that $d_3(t^n\chi)$ represents zero in $E_3$ , which will happen
if either $p$ divides $n$ (in which case $d_4(t^n\chi)=0$) or if 
$\xi'\chi=\xi\chi'$ in $H^*(G)$ for some $\chi'$.  
\par
\proclaim Theorem 3.  In the Lyndon-Hochschild-Serre spectral sequence
with $\fpee$ coefficients for a central extension of a cyclic
$p$-group $C$ by a group $G$, let $u$ generate $E_2^{0,1}$, let $t$
generate $E_2^{0,2}$, and let $d_2(u)=\xi$, $d_3(t)=\xi'$.  \pra
\noindent
a)\qquad If $\chi\in H^*(G)$ is such that $\xi\chi=0$ and $\xi'\chi=0$
in $H^*(G)$, then for all $n$, $t^nu\chi$ survives until $E_4$, and
$$d_4(t^nu\chi)=nt^{n-1}\langle\xi',\chi,\xi\rangle,$$
where $\langle \xi',\chi,\xi\rangle$ is the Massey product (see next
section for a description of the Massey product).
\pra\noindent
b)\qquad If $\chi\in H^*(G)$ is such that $\xi'\chi=\xi\chi'$ for some 
$\xi'$, then for all $n$ $t^n\chi$ survives until $E_4$, and 
$$d_4(t^n\chi)=n(n-1)t^{n-2}u\xi'\chi'.$$
\par\noindent
\remark  The \lq generic' case of statement b) above for $C$ of order
$p$ was conjectured by P.~H.~Kropholler.  He made this
conjecture by looking for a natural map between the relevant
subquotients of $H^*(G)$.  The author thanks P.~H.~Kropholler for
showing him this conjecture, which inspired the work contained in this
paper. 
\par\noindent
\proof  In either case let $\chi$ be an element of $E_0^{*,0}$
yielding the element of the same name in $E_2^{*,0}$, and in case b)
define $\chi'$ similarly.  The proof now splits into two cases.  The
conditions in case a) are equivalent to the existence of $\psi$ and $\psi'$
in $E_0^{*,0}$ such that $d_0(\psi)=0$, $d_0(\psi')=0$,
$d_1(\psi)=\chi\xi$, and $d_1(\psi')=\xi'\chi$.  Define cochains 
$\phi_0,\ldots,\phi_4$ by the equations
\def\sav{(\chi\theta-(-1)^{|\chi|}\psi)}
$$\eqalign{\phi_0&=t^n\chi u,\cr
\phi_1&=\eta_1(n)\chi u+t^n\sav,\cr
\phi_2&=\eta_2(n)\chi u+\eta_1(n)\sav-nt^{n-1}\psi'u,\cr
\phi_3&=\eta_3(n)\chi u+\eta_2(n)\sav
-n\eta_1(n-1)\psi'u-nt^{n-1}\psi'\theta,\cr
\phi_4&=\eta_4(n)\chi u+\eta_3(n)\sav
-n\eta_2(n-1)\psi'u-n\eta_1(n-1)\psi'\theta,}$$
where $\eta_i(j)$ and $\theta$ are as described during and just before
Lemma~2.  The following equations may be seen to hold by applying
Lemma~2. 
$$d_0(\phi_0)=0,\qquad d_1(\phi_i)=d_0(\phi_{i-1}) \quad\hbox{for $i=0,1,2$,}$$
$$d_1(\phi_3)=d_0(\phi_4)-(-1)^{|\chi|}nt^{n-1}(\xi'\psi+\psi'\xi).$$
Now note that in $E_2^{*,*}$ $\phi_0$ represents
$(-1)^{|\chi|}t^nu\chi$, and that $-(\xi'\psi+\psi'\xi)$ represents 
the Massey product $\langle \xi',\chi,\xi\rangle$.  
\par
The proof for case b) is similar.  Here the hypotheses give $\psi\in
E_0^{*,0}$ such that $d_0(\psi)=0$ and $d_1(\psi)=\xi'\chi-\xi\chi'$.  
Now define cochains $\phi_0,\ldots,\phi_4$ by the following equations.
$$\eqalign{\phi_0&=t^n\chi,\cr
\phi_1&=\eta_1(n)\chi-nt^{n-1}u\chi',\cr
\phi_2&=\eta_2(n)\chi-n\eta_1(n-1)u\chi'-nt^{n-1}(\theta\chi'+\psi),\cr
\phi_3&=\eta_3(n)\chi-n\eta_2(n-1)u\chi'-n\eta_1(n-1)(\theta\chi'+\psi),\cr
\phi_4&=\eta_4(n)\chi-n\eta_3(n-1)u\chi'-n\eta_2(n-1)(\theta\chi'+\psi).}$$
Usng Lemma 2 the following equations may be verified.   
$$d_0(\phi_0)=0,\qquad d_1(\phi_i)=d_0(\phi_{i+1}) \quad\hbox{for $i=0,1,2$,}$$
$$d_1(\phi_3)=d_0(\phi_4)-n(n-1)t^{n-2}\xi'u\chi'.$$
\qed
\proclaim Corollary 4.  Let $C$ be a cyclic $p$-group.  Then in the
Serre spectral sequence with $\fpee$ coefficients 
for a fibration with fibre $BC$ and path
connected base space $d_4$ may be described as in Theorem~3.  
\par
\proof For any path connected $X$, Kan-Thurston [\ka] exhibit an aspherical
space $TX$ (in other words a space having $\pi_i(X)=0$ unless $i=1$) 
and a map from $TX$ to $X$ inducing an isomorphism on cohomology for
any local coefficients on $X$.  If we let $X$ be the base space of the 
$BC$ bundle, it follows that the map from $TX$ to $X$ induces an
isomorphism between the Serre spectral sequence for the original
bundle and that for the induced bundle over $TX$.  It follows from the
homotopy long exact sequence that the total space  of the induced
bundle is also aspherical, and so the Serre spectral sequence for this
bundle is isomorphic to the Lyndon-Hochschild-Serre spectral sequence
for the extension given by taking fundamental groups.  
\qed
\par \noindent
\remark The author has not tried to find a direct proof of
Corollary~4.  
\par
\beginsection Massey Products.
\par
We recall first the definition of the Massey (triple) product as in
[\be].  Let $P_\dt$ be any projective $\fpee G$ resolution for $\fpee$, let 
$\Delta$ be a diagonal approximation for $P_\dt$, and let $H$ be a
chain homotopy between $(\Delta\otimes1)\Delta$ and
$(1\otimes\Delta)\Delta$, that is a map from $P_\dt$ to $(P_\dt\otimes
P_\dt\otimes P_\dt)_\dt$ satisfying
$$dH+Hd=(\Delta\otimes1)\Delta-(1\otimes\Delta)\Delta.$$
(Of course such a homotopy exists because $(\Delta\otimes1)\Delta$ and 
$(1\otimes\Delta)\Delta$ are chain maps between projective resolutions 
inducing the same map on homology.)
If $\phi_1$, $\phi_2$ and $\phi_3$ are cochains from $P_\dt$ to
$\fpee$, we may define a triple product $h(\phi_1,\phi_2,\phi_3)$ by
the formula 
$$h(\phi_1,\phi_2,\phi_3)(x)=\mu'(\phi_1\otimes\phi_2\otimes\phi_3)H(x),$$ 
where $\mu'$ is the usual map from $\fpee\otimes\fpee\otimes\fpee$ to
$\fpee$.  If $\phi_1$, $\phi_2$ and $\phi_3$ are cocycles then it is
easily checked that $dh(\phi_1,\phi_2,\phi_3)$ is
$(\phi_1\phi_2)\phi_3-\phi_1(\phi_2\phi_3)$.  If also $\phi_1\phi_2$
and $\phi_2\phi_3$ are coboundaries, with $d\psi_1=\phi_1\phi_2$ and 
$d\psi_2=\phi_2\phi_3$, then it may be verified that the following
cochain is a cocycle.  
$$\langle \phi_1,\phi_2,\phi_3\rangle= (-1)^{|\phi_1|}\phi_1\psi_2
-\psi_1\phi_3+h(\phi_1,\phi_2,\phi_3)$$
For any such $\phi_1$, $\phi_2$ and $\phi_3$ we define the Massey
product of the cohomology classes that they represent to be the
cohomology class of the above cocycle.  If $\Delta$ is strictly
coassociative then $H$ and hence $h(,,)$ may be chosen to be zero.
This is the case in which Massey products are usually discussed,
especially in topology, and is the case used in the previous two
sections of this paper.  From this case it is easy to see that varying
the $\phi_i$ by a coboundary or the $\psi_i$ by a cocycle changes the
Massey product by an element of the ideal of $H*(G)$ generated by
$\phi_1$ and $\phi_3$, so that the Massey product is well-defined
modulo this ideal.  Benson and Evens suggest another extreme case of
interest [\be], the case when $G$ is finite and $P_\dt$ is a minimal
resolution.  In this case the differential on the mod-$p$ cochains on
$P_\dt$ is trivial [\ev], so $\psi_1$ and $\psi_2$ may be chosen to be
zero, and the Massey product is determined by $h$.  We shall use this
case below to determine Massey products in the mod-$p$ cohomology of
finite abelian groups.  First we discuss cyclic groups.  
\par
For $G$ a cyclic group of order $n$ with generator $g$,
Cartan-Eilenberg describe a projective $\fpee G$ resolution $P_\dt$ for
$\fpee$ [\ce] in which $P_i$ is free on one generator $e_i$ with
boundary map 
$$d(e_i)=\cases{(g-1)e_{i-1}&for $i$ odd\cr
\sum_{j=0}^{n-1}g^je_{i-1}&for $i$ even.}$$
Cartan-Eilenberg also give a diagonal approximation for this
resolution whose composite with the projection from $P_\dt\otimes
P_\dt$ to $P_a\otimes P_b$ is given by
$$\Delta (e_{a+b})=\cases{e_a\otimes e_b&for $a$ even,\cr
e_a\otimes ge_b&for $a$ odd, $b$ even,\cr
\sum_{0\leq i<j<n}g^ie_a\otimes g^je_b&for $a$ and $b$ odd.}$$
If $n$ is a power of $p$ then this resolution is minimal.  We now
describe a map $H$ as above for this resolution.  
\par
\proclaim Proposition 5.  Let $P_\dt$ with diagonal approximation
$\Delta$ be the Cartan-Eilenberg resolution for a cyclic group $G$ of
order $n$ generated by $g$ as above.  Define a map $H$ of degree one
from $P_\dt$ to $(P_\dt\otimes P_\dt\otimes P_\dt)_\dt$ by defining
its composite with projection to $P_a\otimes P_b\otimes P_c$ to be 
$$H(e_{a+b+c-1})=\cases{\sum_{0\leq i<j<k<n}g^ie_a\otimes
g^je_b\otimes g^ke_c&if $a$, $b$ and $c$ are all odd,\cr
0&otherwise.}$$
Then $H$ satisfies the formula $dH+Hd=(\Delta\otimes1)\Delta-
(1\otimes\Delta)\Delta$.  
\par
\proof  Check that the maps on each side of the above equation when
composed with projection to $P_a\otimes P_b\otimes P_c$ give the
following map.  
$$e_{a+b+c}\mapsto\cases{
\sum_{0\leq i<j<n}(g^ie_a\otimes g^je_b\otimes e_c-e_a\otimes
g^{i+1}e_b\otimes g^{j+1}e_c) &for $a$, $b$ and $c$ all odd,\cr
\sum_{0\leq i<j<n}g^ie_a\otimes g^je_b\otimes(1-g^{j+1})e_c
&for $a$ and $b$ odd,\cr
\sum_{0\leq i<j<n}g^ie_a\otimes(g^{i+1}-g^j)e_b\otimes g^je_c
&for $a$ and $c$ odd,\cr
\sum_{0\leq i<j<n}(g^i-1)e_a\otimes g^ie_b\otimes g^je_c
&for $b$ and $c$ odd,\cr
0&otherwise.}$$
\qed
\par\noindent
\remark If the group $G$ has order two then $\Delta$ is strictly
coassociative, and the empty sum occuring in the definition of $H$ may
be regarded as zero.  
\proclaim Corollary 6.  Let $G$ be a cyclic group of order $p^m$,
where $p^m > 2$, and let $u$ and $t$ be generators for $H^*(G)$ as in
the introduction, where if $m=1$ then $t$ is the Bockstein of $u$.
The only Massey products that are defined in $H^*(G)$ are $\langle t^iu,
t^ju,t^ku\rangle$, which are defined modulo zero and take the
following value.  
$$\langle t^iu,t^ju,t^ku\rangle= \cases{t^{i+j+k+1}&for $p^m=3$,\cr
0&otherwise.}$$
\par
\proof If $p^m$ is greater than 3 then ${p^m\choose 3}$ is divisible
by $p$, and the result follows easily from Proposition~5.  If $p^m=3$
the argument requires slightly more care (and an explicit definition
of the Bockstein map).  Since this Corollary is well known we omit the
remainder of the proof.  
\qed
\par
\proclaim Proposition 7.  Let $G$ be a finite group expressible as 
$G=G'\times G''$, let $P'_\dt$, $P''_\dt$ be minimal resolutions for
$G'$ and $G''$, and let $h'$, $h''$ be functions as used above to
define Massey products in $H*(G')$ and $H^*(G'')$ respectively.  Then 
$(P'_\dt\otimes P''_\dt)_\dt$ is a minimal resolution for $G$, and the
function $h$ defined (for cochains $\phi_i$ from $P'_\dt$ to $\fpee$
and $\theta_i$ from $P''_\dt$ to $\fpee$) by
$$\eqalign{
h(\phi_1\otimes\theta_1,\phi_2\otimes\theta_2,\phi_3\otimes\theta_3)=
&(-1)^{|\phi_2||\theta_1|+|\phi_3||\theta_1|+|\phi_3||\theta_2|}
(h'(\phi_1,\phi_2,\phi_3)\cr&\otimes(\theta_1\theta_2)\theta_3
+(-1)^{|\phi_1|+|\phi_2|+|\phi_3|}\phi_1(\phi_2\phi_3)\otimes
h''(\theta_1,\theta_2,\theta_3))}$$
may be used to define Massey products in $H^*(G)$.  
\par
\proof  Let $\Delta'$ be a diagonal approximation on $P'_\dt$ and $H'$
a chain homotopy between $(\Delta'\otimes1)\Delta'$ and
$(1\otimes\Delta')\Delta'$ inducing $h'$, and define $\Delta''$ and
$H''$ similarly.  Now the map which is defined to be
$H'\otimes(\Delta''\otimes1)\Delta''+(-1)^i(1\otimes\Delta')\Delta'
\otimes H''$ on $P'_i\otimes P''_j$ is a chain homotopy between 
$(\Delta'\otimes1\otimes\Delta''\otimes1)(\Delta'\otimes\Delta'')$ and 
$(1\otimes\Delta'\otimes1\otimes\Delta'')(\Delta'\otimes\Delta'')$.
The composite of this map with the obvious map from $P_\dt^{\prime\otimes3} 
\otimes P_\dt^{\prime\prime\otimes3}$ to $(P'_\dt\otimes P''_\dt)^{\otimes3}$
which interlaces the factors is therefore a chain homotopy between 
$(\Delta\otimes1)\Delta$ and $(1\otimes\Delta)\Delta$ (where $\Delta 
= (1\otimes\tau\otimes 1)(\Delta'\otimes\Delta'')$), and this is the
map which induces $h$.  
\qed
\par\noindent
\remark  Using Propositions 5 and 7 and Corollary 6 it is easy to
calculate any Massey product in the mod-$p$ cohomology of a finite
abelian group.  We do not include the general formula because it is
long and unilluminating.  Corollary~8 and Theorem~9 contain examples
of such calculations.  
\proclaim Corollary 8.  Let $G$ be a finite abelian group.  Then all
the Massey products that are defined in $H^*(G)$ contain zero unless 
$p=3$ and $G$ has a direct summand of order three.  
\par
\beginsection Examples.  
\par In this short section we consider the 3-groups expressible as an
extension with kernel cyclic of order three and quotient abelian of
rank two.  For these groups we determine all of the differentials in
the Lyndon-Hochschild-Serre spectral sequence with $\fpee$
coefficients, and hence the Poincar\'e series of their cohomology
rings.  The best examples among these (from the point of view of using
our description of $d_4$) are the non-metacyclic examples, but we
shall include the others for completeness.  The
spectral sequence for the nonabelian group of order 27 and exponent 3,
which we shall temporarily refer to as $E$, was first calculated by
Huynh-Mui [\hm].  Our description of $d_4$ greatly simplifies this
calculation.  
The cohomology rings of the metacyclic examples are contained in work of
Diethelm [\di], who used a different spectral sequence.  The other
nontrivial cases of Theorem~9 are new.  
\par
Let $H^*(C_{3^m}\oplus C_{3^n})$ be generated by $y_1$, $y_2$, $x_1$
and $x_2$ of degrees one, one, two and two respectively, where $y_i$
and $x_i$ are in the image of the inflation from projection onto the 
$i$th factor.  The Bockstein maps $y_1$ to $x_1$ (resp.~to zero) if 
$m=1$ (resp. $m>1$), and similarly for $y_2$.  Modulo choice of
generators, we need only consider four distinct extension classes in 
$H^2(C_{3^m}\oplus C_{3^n})$, 0, $x_1$, $x_1+y_1y_2$, and $y_1y_2$.
The first two of these give abelian groups and the third gives a
metacyclic group.  
\par
\proclaim Theorem 9.  Let $G$ be an extension with kernel cyclic of
order three and quotient $C_{3^m}\oplus C_{3^n}$, and let $\xi$ be 
its extension class.  Then the Poincar\'e series for $H^*(G)$ is one
of the following power series.  
$$\halign{\hfil#\qquad&\hfil#\crcr
a) $1/(1-s)^3$,&b) $1/(1-s)^2$,\cr c)
$(1+s)/(1-s)(1-s^6)$,& d) $(1+s+s^2)/(1-s^2)^2(1-s)$,\cr
e) $(1+s^2)/(1-s^6)(1-s)^2$,&f)
$(1+s+2s^2+2s^3+s^4+s^5)/(1-s^6)(1-s)$.\cr}$$
If $\xi=0$ then case a) occurs.  If $\xi=x_1$ then case b) occurs.  
If $\xi=x_1+y_1y_2$ then case b) occurs if $n$ is greater than one and 
case c) occurs if $n=1$.  If $\xi=y_1y_2$ then d) occurs if both $m$
and $n$ are greater than one, e) if one of $m$ and $n$ is one, and f)
if $m$ and $n$ are both one.  
\par
\noindent
\proof  As stated above, we examine the Lyndon-Hochschild-Serre
spectral sequences for these extensions.  It is convenient to
introduce the bigraded Poincar\'e series for the pages of the spectral
sequence, defined by the following formula.
$$P_r(s,s')=\sum_{i,j}s^is'^j\dim E_r^{i,j}$$
The Poincar\'e series that we wish to determine is equal to
$P_\infty(s,s)$.  For any central extension of $p$-groups the image of
$d_r$ is a finitely generated graded module for the kernel of $d_r$
which in turn is a finitely generated graded $\fpee$-algebra.  It
follows that $P_r$ and $P_{r+1}$ satisfy a relation of the following form,
$$P_r(s,s')-P_{r+1}(s,s')=(s^r+s'^{r-1})f(s,s')/g(s,s'),$$
where $f(s,s')$ is an integer polynomial and $g(s,s')$ is a product of
terms of the forms $1-s^k$ and $1-s'^l$.  This relation provides a
useful check when calculating $P_r(s,s')$.  
\par
The $E_2$-page of the spectral sequence we are considering is 
$\Lambda[u,y_1,y_2]\otimes\fpee[t,x_1,x_2]$, and
$P_2=1/(1-s)^2(1-s')$.  All of the groups have an abelian subgroup of
index three and an element in the cohomology of this subgroup that
maps to $t$ in the cohomology of the kernel.  A simple Evens' norm map
argument [\ev] now implies that in the spectral sequence $t^3$ is a
permanent cycle, and so $E_7= E_\infty$.  The first two differentials
are given by $d_2(u)=\xi$, $d_2(t)=0$, $d_3(t)=\beta(\xi)$.  
\par
If $\xi$ is zero then the spectral sequence collapses.  If $\xi$ is
either $x_1$ or $x_1+y_1y_2$ then $d_2$ is injective from odd rows to
even rows, and $P_3=(1+s)/(1-s'^2)(1-s)$.  Now $d_3(t)$ is zero if 
either
$\xi$ is $x_1$ or $\xi$ is $x_1+y_1y_2$ and $m$, $n$ are both greater
than one.  If $m=1$ and $n>1$ then
$\beta(x_1+y_1y_2)=(x_1+y_1y_2)y_2$, so that $d_3(t)$ is zero in this
case too.  In the other cases when $\xi=x_1+y_1y_2$, $d_3(t)$ is equal
to $-x_2y_1$ modulo the image of $d_2$, and in these cases 
$P_4=(1+(s-s^3)(1+s'^2)+ss'^4)/(1-s'^6)(1-s)$.  The $E_4$-page is
concentrated in even vertical degree, so $d_4$ must be zero.  The only
possible non-zero differential now is $d_5(t^2y_1)$, but
$d_5(t^2x_2y_1)$ is zero by Kudo's transgression theorem [\ku], and
hence $d_5(t^2y_1) $ is zero since multiplication by $x_2$ is
injective from $E_5^{6,0}$ to $E_5^{8,0}$.  Hence $E_4=E_\infty$.  
\par
We now move on to the cases when $\xi=y_1y_2$.  Here
$P_3$ is equal to $(1+2s+s'(2s+s^2))/(1-s^2)(1-s^2)^2$.   
If $m$ and $n$ are both greater than one, then $d_3(t)$ is zero and so
$E_3=E_\infty$.  If exactly one of $m$ and $n$ is equal to one, we may
assume that $n=1$.  Now $d_3(t)= -x_2y_1$, and 
$$\eqalign{P_4=
(1+2s-s^3&+s'(2s+s^2-s^4)+s'^2(2s-s^3)\cr
&+s'^3(s+s^2-s^4)+2s'^4s+s'^5(s+s^2))/(1-s'^6)(1-s^2)^2.\cr}$$
We now need to determine the images of $tuy_1$, $t^2y_1$, $t^2y_2$ and
$t^2uy_1$ under $d_4$.  From Proposition~7 and Corollary~6 it follows
that 
$$\langle x_2y_1,y_1,y_1y_2\rangle = \langle y_1,y_1,y_1\rangle x_2y_2=0,$$
and so by Theorem~3, $d_4(t^iuy_1)$ is zero.  Also $(x_2y_1)y_1=0$,
whereas $(-x_2y_1)y_2=-(y_1y_2)x_2$, and so Theorem~3 implies that
$d_4(t^2y_1)=0$ and $d_4(t^2y_2)=-ux_2^2y_1$.  The Poincar\'e series
for the $E_5$-page is given by the following formula.  
$$\eqalign{P_5=
(1+2s-s^3&+s'(2s+s^2-s^4-s^5)+s'^2(2s-s^3)\cr&
+s'^3(s+s^2-s^4)+s'^4s+s'^5(s+s^2))/(1-s'^6)(1-s^2)^2.\cr}$$
The remaining possibly non-zero differentials are $d_5(t^2y_1)$ and
$d_5(t^2uy_1)$, which we claim are zero, and then $d_6(t^2uy_1)$,
which we also claim to be zero.  It may be shown that $d_5(t^2y_1)=0$
either using a Kudo transgression theorem argument or by noting that
no polynomial in $x$ and $x'$ can be hit, since $E_5^{*,j}=E_6^{*,j}$
for $j$ equal to 2 and 3, and each of these rows contains a free
submodule for the subring of $E_5^{*,0}$ generated by $x$ and $x'$.  
Showing that $d_5(t^2uy_1)$ is zero is more difficult.  Since
$E_5^{6,1}$ is generated by $uy_1y_2x_2^2$, it suffices to show that a
cohomology element may be chosen which yields $uy_1y_2x_2^2$ and whose
Bockstein yields a non-zero element in $E_\infty^{7,1}$ (which is
clearly isomorphic to $E_5^{7,1}$).  Let $\hat E_r^{i,j}$ stand for
the corresponding spectral sequence with integer coefficients.  There
is a Bockstein map from $\hat E_r^{i,j}$ to $\hat E_r^{i,j+1}$ if
$j>0$ (resp.~$\hat E_r^{i+1,j}$ if $j=0$) converging to the usual
Bockstein from mod-$p$ cohomology to integral cohomology.  Under this
Bockstein $uy_1y_2$ maps to a generator for $\hat E_\infty^{2,2}$ (one
need only compute $d_3$ to check this), which is clearly not divisible
by three as a cohomology element because $\hat E_\infty^{1,3}$ and 
$\hat E_\infty^{0,4}$ are trivial.  It follows that the Bockstein of
an element yielding $uy_1y_2$ must yield a nonzero multiple of 
$ux_2y_1$ in $E_*$, and hence that $d_5(t^2uy_1)$ cannot be
$uy_1y_2x_1^2$, so must be zero.  Similarly, $d_6(t^2uy_1)$ must be in
the kernel of the Bockstein from $E_6^{7,0}$ to $E_6^{8,0}$, so can
only be a multiple of $x_1^3y_1$.  Comparison with the spectral
sequence for the subextension with quotient $C_{3^m}\oplus\{0\}$ shows
now that $d_6(t^2uy_1)=0$.  
\par
Now we consider the case when $\xi=y_1y_2$ and $m=n=1$.  Here $d_3(t)$
is equal to $x_1y_2-x_2y_1$, and $P_4$ is given by the following formula.
$$\eqalign{
P_4=(1+2s-s^3&+s'(2s+s^2-2s^4+s^6)+s'^2(2s-s^3)\cr&+s'^3(s^2+s^3-2s^4+s^6)
+2s'^4s+s'^5(s^2+s^3))/(1-s^2)^2(1-s'^6)\cr}$$
We need to determine the images under $d_4$ of $tu(x_1y_2-x_2y_1)$, 
$t^2y_1$, $t^2y_2$ and $t^2u(x_1y_2-x_2y_1)$.  From Propostion~7 and
Corollary~6 we see that 
$$\langle x_1y_2-y_1x_2,x_1y_2-y_1x_2,y_1y_2\rangle=
\langle y_1,y_1,y_1\rangle x_2^2y_2 - x_1^2y_1\langle
 y_2,y_2,y_2\rangle = x_1x_2^2y_2-x_1^2x_2y_1,$$
and hence that
$d_4(t^iu(x_1y_2-x_2y_1))=it^{i-1}(x_1x_2^2y_2-x_1^2x_2y_1)$.  Note
that modulo multiples of $x_1y_2-x_2y_1$, $x_1x_2^2y_2-x_1^2x_2y_1$ is
equal to the first reduced power of $x_1y_2-x_2y_1$, which we know
must be hit by some differential.  Also note that 
$$(x_1y_2-x_2y_1)y_1=-(y_1y_2)x_1\qquad(x_1y_2-x_2y_1)y_2=-(y_1y_2)x_2,$$
from which it follows that $d_4(t^2y_i)=u(x_1y_2-x_2y_1)x_i$.  
The Poincar\'e series for $E_5$ page is the following expression. 
$$\eqalign{
P_5=(1+2s-s^3-s^7&+s'(2s+s^2-2s^4-2s^5+s^6)+s'^2(2s-s^3-s^7)
\cr &+s'^3(s^2-2s^4+s^6)
+s'^4s+s'^5s^2)/(1-s^2)^2(1-s'^6)\cr}$$
It remains now to calculate $d_5(t^2(x_1y_2-x_2y_1))$, which is
$x_1^3x_2-x_2^3x_1$ by Kudo's transgression theorem, and
$d_5(t^2uy_1y_2)$, which must be a non-zero multiple of
$u(x_1^3y_2-x_2^3y_1)$ to ensure that $E_6^{*,1}$ has a well defined
module structure for $E_6^{*,0}$.  The Poincar\'e series for $E_6$,
which is equal to $E_\infty$, is the following. 
$$\eqalign{
P_6=(1+2s+s^2&+s^3+s^4+s^5+s^6+s'(2s+s^2+2s^3-s^4)\cr
&+s'^2(2s+s^3+s^5)
+s'^3(s^2-s^4))/(1-s^2)(1-s'^6)\cr}$$
\qed
\par

\beginsection References. 

\frenchspacing
\item{[\be]} D. J. Benson, 
{\it Representations and Cohomology I}, Cambridge Studies in
Advanced Mathematics, Cambridge University Press (1991).
\par
\item{[\ca]} H. Cartan,
Sur les groupes d'Eilenberg-MacLane $H(\pi,n)$ I and II, 
{\it Proc. Nat. Acad. Sci. Amer.}  {\bf 40} (1954), 467--471 and 704--707.
\par
\item{[\ce]} H. Cartan and S. Eilenberg,
{\it Homological Algebra}, Princeton University Press (1956).
\par
\item{[\di]} T. Diethelm, The mod-$p$ cohomology rings of the non
abelian split metacyclic $p$-groups, {\it Arch. Math. Basel} {\bf 44}
(1985) 29--38.
\par
\item{[\ev]} L. Evens, {\it The Cohomology of Groups}, Oxford
Mathematical Monographs, Oxford University Press (1991).
\par
\item{[\hm]} Hu\`ynh M\`ui, The mod-$p$ cohomology algebra of the
extra-special group $E(p^3)$, preprint (1983).  
\par
\item{[\ka]} D. M. Kan and W. P. Thurston, Every connected space has the
homology of a $K(\pi,1)$, {\it Topology} {\bf 15} (1976) 253--258.
\par
\item{[\ku]} T. Kudo, A transgression theorem, {\it Mem. Fac. Sci.
Kyusyu Univ. Series A}, {\bf 9} (1956) 79--81.
\par
\item{[\se]} J. P. Serre,
Cohomologie modulo 2 des complexes d'Eilenberg-MacLane, 
{\it Comm. Math. Helv.} {\bf 27} (1953) 192--232.
\par
\item{[\st]} N. E. Steenrod, Products of cocycles and extensions of
mappings, {\it Ann. Math.} {\bf 48} (1947) 290--320.
\par

\bye